\documentclass [12pt,leqno] {article}

\usepackage{amsfonts,url}

\def\sgn{\mathop{\rm sgn}\nolimits}
\def\dist{\mathop{\rm dist}\nolimits}
\def\supp{\mathop{\rm supp}\nolimits}

\newcommand{\qed}{~\hfill~$\fbox{}$ }

\newcommand{\proces}{\{ X_t:\:t\geq 0\}} 
\newcommand{\RR}{\mathbb{R}}
\newcommand{\R}{ \mathbb{R}^{d}}
\newcommand{\N}{\mathbb{N}}

\newcommand{\dowod}{{\em Proof}.\/ }
\newcommand{\indyk}[1]{{\bf 1}_{#1}}
\newcommand{\sfera}{ \mathbb{S}}

\newcommand{\Fourier}{ {\cal F}}

\newcommand{\nubounded}[1]{\bar{\nu}_{#1}}

\newcommand{\scalp}[2]{#1\cdot#2}

\newcommand{\numj}{\nu}

 \def\dist{\mathop{\rm
    dist}\nolimits} \def\diam{\mathop{\rm diam}\nolimits}

\newtheorem{lemat}{\indent\sc Lemma}

\newtheorem{twierdzenie}{\indent\sc Theorem}
\newtheorem{wniosek}[lemat]{\indent\sc Corollary}

\newcounter{conum} 

\begin{document}

\title{Transition density estimates for jump L\'evy processes}
\author{Pawe{\l} Sztonyk}
\footnotetext{ Institute of Mathematics and Computer Science,
  Wroc{\l}aw University of Technology,
  Wybrze{\.z}e Wyspia{\'n}\-skie\-go 27,
  50-370 Wroc{\l}aw, Poland\\
  {\rm e-mail: sztonyk@pwr.wroc.pl} \\
}
\date{June 29, 2010}
\maketitle

\begin{center}
  Abstract
\end{center}
\begin{scriptsize}
  Upper estimates of densities of convolution semigroups of probability measures are given under explicit assumptions
  on the corresponding L\'evy measure and the L\'evy--Khinchin exponent.
\end{scriptsize}

\footnotetext{2000 {\it MS Classification}:
Primary 60G51, 60E07; Secondary 60J35, 47D03, 60J45 .
The research was supported by grant MNiSW N N201 397137.
}
\section{Introduction}\label{Intro}

We investigate a class of jump--type L\'evy processes and the corresponding convolution semigroups of measures and give estimates of 
their transition densities. 
The primary examples of considered objects are the stable L\'evy processes. Stable L\'evy processes are important in theoretical probability, physics and finance, and their asymptotic properties are of interest for many authors. The estimates of their transition densities were obtained, e.g., in \cite{BlGe,PT,Hi1,Hi2,G,GH,W,BS2007}. Our present study applies to the more general layered and tempered stable processes 
(for definitions see Section \ref{sec:examples}) with marginal tails
heavier than Gaussian but lighter than stable processes. The layered and tempered L\'evy processes 
are known in physical literature as {\it truncated L\'evy flights} 
and were introduced in statistical physics to model turbulence (\cite{Novikov,Mantegna,Koponen}). Numerical analysis of stock price time series also indicate distributions of this type (\cite{WTT}). In particular they are used
to model stochastic volatility (\cite{CGMY1,CGMY2}). Our main objective is to describe the radial decay of the transition density when the radial decay {\it and} some degree of smoothness of the L\'evy measure are given. In particular we obtain precise off-diagonal estimates of the transition density even if the L\'evy measure is singular. We will now describe our results.

Let $d\in\{1,2,\dots\}$, $b\in\R$, and $\nu$ be a L\'evy measure on $\R$, i.e.,
$$
  \int_{\R} \left(1\wedge |y|^2\right)\,\nu(dy) < \infty.
$$ 

We consider the convolution semigroup of probability measures $\{ P_t, t\geq 0 \}$ with the Fourier transform $\Fourier(P_t)(\xi)=\exp(-t\Phi(\xi))$, where
\begin{equation}\label{eq:Phi_2}
  \Phi(\xi) =   - \int \left(e^{i\scalp{\xi}{y}}-1-i\scalp{\xi}{y}\indyk{B(0,1)}(y)\right)\nu(dy) - i\scalp{\xi}{b}  ,\quad \xi\in\R.
\end{equation}
In Section \ref{sec:construction} we 
recall a procedure of the construction of the semigroup. The main result of the present paper is the following theorem.

\begin{twierdzenie}\label{th:main} Assume that
$\nu$ is a L\'evy measure such that
\begin{equation}\label{eq:nu_intro}
  \nu(A) \leq \int_\sfera\int_0^\infty \indyk{A}(s\theta) s^{-1-\alpha} q(s)\phi(s) \, ds\mu(d\theta),\quad A\in{\cal B}(\R ),
\end{equation}
where $\alpha\in(0,2)$, $\mu$ is a finite nonnegative measure on the unit sphere $\sfera$, 
and $q:\:[0,\infty)\to (0,\infty)$ is bounded nonincreasing function such that
\begin{equation}\label{eq:q_intro}
  q(s) \leq \kappa_1 q(2s),\quad s>0,
\end{equation}
for a positive constant $\kappa_1$,
and 
$\phi:\:[0,\infty)\to (0,\infty)$
is a bounded nonincreasing function satisfying
\begin{equation}\label{eq:phi_intro}
  \phi(a)\,\phi(b) \leq  \kappa_2 \phi(a+b),\quad a,b \geq 0,
\end{equation}
where $\kappa_2>0$.

Furthermore we assume that there exists $\beta\in [\alpha,2]$ and finite $c\geq 0$ such that
\begin{equation}\label{eq:beta_intro}
  \int_0^r s^{1-\alpha} q(s)\frac{\phi(s)}{\phi(s/2)}\, ds \leq c r^{2-\beta},\quad r>1,
\end{equation}
and
\begin{equation}\label{eq:Ass0}
  \Re(\Phi(\xi))\geq c \left( |\xi|^{\alpha}\wedge |\xi|^\beta \right),\quad \xi\in\R.
\end{equation}

We finally assume that $\mu$ is a $\gamma-1$-measure on $\sfera$ for some $\gamma\in [1,d]$, i.e.
\begin{equation}\label{eq:mu_gamma_measure}
  \mu(\sfera\cap B(\theta,\rho)) \leq c \rho^{\gamma -1},\quad \theta\in\sfera,\,\rho>0.
\end{equation}

Then the measures $P_t$ are absolutely continuous with respect to the Lebesgue measure and their densities $p_t$ satisfy
\begin{eqnarray*}
  p_t(x+tb_{t^{1/\alpha}}) 
  & \leq & c t^{-d/\alpha} \min\left\{1, t^{1+\frac{\gamma}{\alpha}}
          |x|^{-\gamma-\alpha}q(|x|)\phi\left(|x|/4\right) 
            \right\},
\end{eqnarray*}
for $t\in(0,1],\,x\in\R,$
and
\begin{eqnarray}\label{eq:es_above}
  p_t(x+tb_{t^{1/\beta}}) 
  & \leq & c_1 t^{-d/\beta} \left(\min\left\{1, t^{1+\frac{\gamma}{\beta}}
          |x|^{-\gamma-\alpha}q(|x|)\phi\left(|x|/4\right)\right\} \right.\nonumber \\ 
  &      & \left. +\, e^{-c_2t^{-1/\beta}|x|\log(1+c_3t^{-1/\beta}|x|)} \right),
\end{eqnarray}
for $t> 1,\,x\in\R$, where
\begin{equation}\label{eq:def_br}
  b_r = \left\{
  \begin{array}{ccc}
    b - \int_{r<|y|<1} y \, \nu(dy) & \mbox{  if  } & r \leq 1,\\
    b + \int_{1<|y|<r} y \, \nu(dy) & \mbox{  if  } & r > 1.
  \end{array}\right.
\end{equation}
\end{twierdzenie}

In polar coordinates we have $\nu(ds d\theta)\leq s^{-1-\alpha} q(s)\phi(s)ds\mu(d\theta)$. The {\it relative} radial decay of the upper bound is the same  in each direction, but its {\it absolute} size is highly anisotropic, as permitted by the spectral measure $\mu$. Condition (\ref{eq:Ass0}) is our only lower bound for $\nu$. In particular, if $\mu$ is nondegenerate, i.e., the support of $\mu$ is not contained in any proper linear subspace
of $\R$, and for a constant $r_0>0$ we have
$$
  \nu(ds d\theta)\geq c \indyk{(0,r_0]}(s)s^{-1-\alpha} ds\mu(d\theta),
$$ then $\Re(\Phi(\xi))\geq c \left( |\xi|^{\alpha}\wedge |\xi|^2 \right)$. Furthermore, if for $\beta\in [\alpha,2)$ we have
$$
  \nu(ds d\theta)\geq c \left[\indyk{(0,r_0]}(s)s^{-1-\alpha} ds + \indyk{(r_0,\infty)}(s)s^{-1-\beta} \right]\mu(d\theta),
$$
then $\Re(\Phi(\xi))\geq c \left( |\xi|^{\alpha}\wedge |\xi|^\beta \right)$.
  
We should also note that the {\it doubling property} (\ref{eq:q_intro}) is equivalent to the following:
there exist $c>0$ and $\eta\geq 0$ such that
\begin{equation}\label{eq:q_intro_2}
   \frac{q(r)}{q(R)} \leq c\left(\frac{r}{R}\right)^{-\eta},\quad 0<r\leq R.
\end{equation}
The typical examples of such functions $q$ are $q(s)=(1+s)^{-a},$ or $q(s)=(\log(e+s))^{a}(1+s)^{-ma},$ for $a\geq 0$ and
$m> 1$.

If a function $\phi$ satisfies (\ref{eq:phi_intro}) then $1/\phi$ is submultiplicative in the sense of Section 25 in \cite{Sato} and locally bounded.
By \cite[Lemma 25.5]{Sato} there exist $c_1,c_2$ such that $\phi(s)\geq c_1 e^{-c_2s}$, $s\geq 0$. Functions which satisfy
(\ref{eq:phi_intro}) are, e.g., $\exp(-ms^{a})$, $1/\log^m(s+e)$
for every $m>0$ and $0<a \leq 1$. A product of two functions which satisfy (\ref{eq:phi_intro}) also satisfies (\ref{eq:phi_intro}) 
(see \cite[Proposition 25.4]{Sato}).

We note that (\ref{eq:mu_gamma_measure}) is equivalent to $\nu$ being a $\gamma$-measure near $S$: $\nu(B(\theta,r))\leq c r^{\gamma}$ 
for $\theta\in S$, $r<1/2$.

We would like to mention recent results related to Theorem \ref{th:main}. The rotation invariant $\alpha$-stable L\'evy process has the L\'evy measure $\nu(dy)=c|y|^{-d-\alpha}$. The asymptotic behaviour of its densities
is well known (see, e.g., \cite{BlGe}) and in this case we have $p_t(x)\approx \min(t^{-d/\alpha},t|x|^{-d-\alpha})$. W.E. Pruitt and S.J. Taylor investigated in \cite{PT} multivariate stable densities in a more general setting. They obtained the upper bound 
$p_t(x)\leq ct^{-d/\alpha} (1+t^{-1/\alpha}|x|)^{-1-\alpha}$ by Fourier-analytic methods. We like to note that such decay is indeed observed if the spectral measure $\mu$ of the stable L\'evy has an atom (for lower bounds see \cite{Hi1,Hi2} and \cite{W}). Using the perturbation formula P.~G{\l}owacki and W. Hebisch proved in \cite{G} and \cite{GH} that if 
$\nu(dx) = |x|^{-\alpha-d}g(x/|x|)dx$ and $g$ is bounded, then $p_t(x)\leq c \min(t^{-d/\alpha},t|x|^{-d-\alpha})$. 
When $g$ is continuous on $\sfera$ we even have
$\lim_{r\to\infty}r^{d+\alpha}p_1(r\theta)=cg(\theta)$, $\theta\in\sfera$ and if $g(\theta)=0$ then additionally
$\lim_{r\to\infty}r^{d+2\alpha}p_1(r\theta)=c_{\theta}>0$, which was proved by J. Dziuba\'nski in \cite{Dziub}.
A. Zaigraev in \cite{Zai} obtained further asymptotic expansions of the $\alpha$--stable density for sufficiently regular $g$. K.~Bogdan and T.~Jakubowski in \cite{BJ2007} obtained estimates of heat kernels of the fractional Laplacian perturbed by gradient operators. Derivatives of 
stable densities have been considered in \cite{Lewand} and \cite{Szt2007}. The existence and upper estimates of the densities for symmetric L\'evy and L\'evy--type
processes were investigated also by V. Knopova and R.L. Schilling in \cite{KnopSch1,KnopSch2}.

More recent asymptotic results for stable L\'evy processes are given 
in papers \cite{W} and \cite{BS2007}. In particular if the L\'evy measure is given in polar coordinates as
$\nu(dsd\theta)=s^{-1-\alpha}ds\mu(d\theta)$ and $\mu$ is symmetric and satisfies (\ref{eq:mu_gamma_measure}) for some $\gamma\in[1,d]$ 
then we have $p_t(x)\leq c \min\{t^{-d/\alpha},t^{1+(\gamma-d)/\alpha}|x|^{-\gamma-\alpha}\}$.

In the present paper we strengthen the results and the methods of
\cite{Sztonyk2010_2}, which was restricted to symmetric processes and $\phi\equiv 1$. 
We take this opportunity to notice that one of our estimates in \cite{Sztonyk2010_2} is incorrect. Namely, the upper bound of the density
for large times $t$ should contain $t^{1+\frac{\eta+\gamma-d}{\beta}}$ instead of
$t^{1+\frac{\gamma-d}{\beta}}$ according to the proof given in \cite[Theorem 3]{Sztonyk2010_2} but the conclusion of the proof failed to summarize the estimates properly. Here we override the estimates with the more adequate (\ref{eq:es_above}).

We also want to mention estimates for the transition density of a class of Markov processes with jump intensities which are not necessarily translation invariant but dominated
by the L\'evy measure of the stable rotation invariant process given in \cite{Sztonyk2010_3} (see also \cite{ShUe}). It is an open and important problem to extend the present techniques to the Markov (not translation-invariant) case. 
Existing results for Markov transition semigroups (\cite{ChKum,ChKum08}) are mostly based on assuming that the L\'evy (jump) kernel of 
the Markov process has isotropic upper and lower bounds in small scales. 
In sharp contrast, our assumptions allow for a highly anisotropic behavior of the L\'evy measure $\nu$. 
This is so because the spectral measure $\mu$ in (\ref{eq:nu_intro}) does not in general have a finite isotropic majorant, 
unless $\mu$ is (bounded by) the uniform measure on the unit sphere, which is a trivial case in our considerations.

The paper is organized as follows. In Section \ref{sec:construction} we give preliminaries.
Section \ref{sc:Trunc} contains estimates of the transition densities for the L\'evy measures with bounded support. 
In Section \ref{sc:Large_Jumps} we investigate asymptotics of the convolution exponents of finite measures. In Section \ref{sc:Main_Theorem} we
prove Theorem \ref{th:main} by combining the results of Section 4 and Section 5. In Section \ref{sec:examples} we give examples of processes satisfying the assumptions of Theorem \ref{th:main}.


\section{Construction of the semigroup}\label{sec:construction}

For $x\in\R$ and $r>0$ we let $|x|=\sqrt{\sum_{i=1}^d x_i^2}$ and 
$B(x,r)=\{ y\in\R :\: |y-x|<r \}$. We denote $\sfera=\{ x\in\R:\: |x|=1 \}$.
All the sets, functions and measures considered in the sequel will be Borel.
For a measure $\lambda$ on $\R$, $|\lambda|$ denotes
its total mass.  For a function $f$ we let $\lambda(f)=\int f d\lambda$,
whenever the integral makes sense.
When $|\lambda|<\infty$ and $n=1,2,\ldots$ we let $\lambda^{n*}$
denote the $n$-fold convolution of $\lambda$ with itself:
$$
  \lambda^{n*}(f)=\int f(x_1+x_2+\dots+x_n)\,\lambda(dx_1)\lambda(dx_2)\ldots\lambda(dx_n)\,.
$$
We also let $\lambda^{0*}=\delta_0$, the evaluation at $0$, and
$$
  e^{\lambda}(f)=\sum_{n=0}^\infty \frac{\lambda^{n*}(f)}{n!}.
$$

Let $\nu$ be a L\'evy measure on $\R$, i.e., $\int \left(|y|^2\wedge 1\right) \nu(dy)<\infty$. For all $0\leq \rho < r \leq\infty$ we denote 
$\numj_{\rho,r}(dy) = \indyk{B(0,r)\setminus B(0,\rho)}(y)\,\nu(dy)$.

We now construct the semigroup $\{ P_t,\,t\geq 0\}$ corresponding to $\nu$ (for a more axiomatic
introduction to convolution semigroups we refer the reader to \cite{BgFt, Jc1}). We denote
$$
  a_{\rho,r}= - \int_{B(0,1)} y \,\nu_{\rho,r}(dy),\quad 0 < \rho < r \leq\infty,
$$
and consider the probability
measures
$$
  P_t^{\rho,r} = e^{t(\numj_{\rho,r}-|\numj_{\rho,r}|\delta_0)} * \delta_{ta_{\rho,r}},\quad 0 < \rho < r \leq\infty,\,t\geq 0.
$$
$P^{\rho,r}_t$ form a convolution semigroup:
$$
  P^{\rho,r}_t*P^{\rho,r}_s=P^{\rho,r}_{s+t}\,,\quad s,\,t\geq 0\,.
$$
The Fourier transform of $P^{\rho,r}_t$ is
$$
  \Fourier(P^{\rho,r}_t)(\xi)=
  \int e^{ i\scalp{\xi}{y} } P^{\rho,r}_t (dy) =
  \exp \left(t\int \left(e^{ i\scalp{\xi}{y} }-1-i\scalp{\xi}{y} \indyk{B(0,1)}(y)\right)\,\numj_{\rho,r}(dy)\right),
  \quad \xi\in\R\,.
$$
We have
\begin{equation}\label{eq:var}
  \int |y|^2 \,P^{\rho,r}_t(dy) = \sum_{k=1}^{d} \left[-\frac{\partial^2}{\partial \xi_k^2}\Fourier(P^{\rho,r}_t)(\xi)\right]_{\xi=0}=t \int |y|^2\,\numj_{\rho,r}(dy),
\end{equation}
for all $0< \rho < r \leq 1$. Furthermore $P^{\rho,1}_t=P^{r,1}_t * P^{\rho,r}_t$, and this yields
\begin{eqnarray*}
  P^{\rho,1}_t - P^{r,1}_t 
  &   =  & P^{r,1}_t * P^{\rho,r}_t  - P^{r,1}_t \\
  &   =  & P^{r,1}_t * \left(P^{\rho,r}_t  - \delta_0 \right).
\end{eqnarray*}
Let $f$ be a continuous bounded function on $\R$.  
For every $\varepsilon>0$ and $R>0$ we can choose $\delta>0$ such that $|f(x+y)-f(x)|<\varepsilon/3$ for all $x\in B(0,R)$, $|y|<\delta$,  and by 
(\ref{eq:var}) and Chebyshev's inequality we get
\begin{eqnarray*}
  \left|\int f(x+y) \,P^{\rho,r}_t (dy) - f(x)\right| 
  &   =  & \left|\int \left(f(x+y)- f(x)\right)\, P^{\rho,r}_t (dy)\right| \\
  & \leq & \int \left|f(x+y)-f(x)\right|\, P^{\rho,r}_t (dy) \\
  & \leq & \int_{|y|<\delta} \frac{\varepsilon}{3} \, P^{\rho,r}_t (dy) + \int_{|y|\geq\delta} 2\|f\|_\infty\, P^{\rho,r}_t (dy) \\
  & \leq & \frac{\varepsilon}{3} + 2\|f\|_\infty  \frac{t \int |y|^2\,\numj_{0,r}(dy)}{\delta^2}.
\end{eqnarray*}

Taking $R>0$ such that $R^2> \frac{6\|f\|_\infty t\int |y|^2\,\numj_{0,1}(dy)}{\varepsilon}$ and $r>0$ such that $\int |y|^2\,\numj_{0,r}(dy)<\frac{\delta^2}{2t\|f\|_\infty}\frac{\varepsilon}{3}$ we obtain
\begin{eqnarray*}
  \left|\left(P^{\rho,1}_t - P^{r,1}_t\right)(f)\right|
  & \leq & \int \left(\int \left|f(x+y)-f(x)\right| \,P^{\rho,r}_t(dy)\right) \,P^{r,1}_t(dx) \\
  & \leq & \int_{B(0,R)}\left(\int \left|f(x+y)-f(x)\right| \,P^{\rho,r}_t(dy)\right)\,P^{r,1}_t(dx) \\
  &      & +\int_{B(0,R)^c}2\|f\|_\infty \,P^{r,1}_t(dx) \\
  & \leq & \frac{\varepsilon}{3} + 2\|f\|_\infty  \frac{t \int |y|^2\,\numj_{0,r}(dy)}{\delta^2} 
           + 2\|f\|_\infty \frac{t \int |y|^2\,\numj_{0,1}(dy)}{R^2} \leq \varepsilon.
\end{eqnarray*}

It follows that the measures $P^{r,1}_t$ weakly
converge to a probability measure $P^{0,1}_t$ as $r \to 0$.
We let $P_t=P^{0,1}_t*P^{1,\infty}_t*\delta_{tb}$.
$\{ P_t,\, t\geq 0 \}$ is also a convolution semigroup and
$\Fourier(P_t)(u)=\exp(-t\Phi(u))$, where
$$
  \Phi(\xi)  =   -\int \left(e^{ i\scalp{\xi}{y} }-1-i\scalp{\xi}{y}\indyk{B(0,1)}(y)\right)\nu(dy) - i\scalp{\xi}{b},\quad \xi\in\R.
$$
We call $\nu$ the {\it L\' evy measure} of the semigroup $\{P_t\,,\;t\geq
0\}$ \cite{Ho, BgFt}. The semigroup determines the stochastic L\'evy process $\proces$ on $\R$ with transition probabilities $P(t,x,A)=P_t(A-x)$.

In what follows {\it constants} mean positive real numbers depending only on $d,\alpha,\beta,\mu,q,\phi,\gamma,\eta$.  
We write $f\approx g$ to indicate that there is a constant 
$c$ such that $c^{-1}f \leq g \leq c f$.

\section{L\'evy measure with bounded support}\label{sc:Trunc}

We consider in this section a L\'evy measure $\nu_0$ such that $\supp \nu_0\subset B(0,r)$ for some $r>0$. Let $\pi$ be a infinitely divisible
distribution with generating triplet (see \cite{Sato}) $(0,\nu_0,0)$, i.e.
 
\begin{equation}\label{eq:F_pi}
  \Fourier(\pi)(\xi) =
  \exp\left( \int (e^{i\scalp{\xi}{y}}-1-i\scalp{\xi}{y}\indyk{B(0,1)}(y))
  \,\nu_0(dy)\right)\, ,
  \quad \xi\in\R\, .
\end{equation}

The L\'evy measures with bounded support are discussed, e.g., in Section 26 of \cite{Sato}. It follows from \cite[Theorem 26.1]{Sato} that
$\pi(B(0,a)^c) = o(e^{-\kappa a\log a})$ as $a\to\infty$ for all $\kappa\in(0,1/r)$.
In this section we will essentially use the methods of \cite{Sato} but our 
computations have to be a little bit more precise. 
In particular, as a consequence of the estimate of tails given in Lemma \ref{lm:pi_tail} below we obtain in Lemma \ref{lm:pi_estimate} an estimate of a density of $\pi$. 
We note that sharp estimates for the transition densities of the truncated stable processes 
(with $\nu(dx)=c|x|^{-d-\alpha}\indyk{|x|<1}dx$) were given in \cite{ChKK}.

We let
$$
  \xi_0^j=\int_{|y|>1} y_j\, \nu_0(dy),\quad j=1,2,...,d, 
$$
$$
  \xi_0=\max\{ |\xi_0^1|,\dots,|\xi_0^d| \},
$$
and
$$
  M_j=\int y_j^2\, \nu_0(dy), \quad j=1,2,...,d,
$$
$$
  M=\max\{M_1,\dots,M_d\}.
$$

\begin{lemat}\label{lm:pi_tail}
  For every $a\geq 2\sqrt{d}\left(\xi_0+M/e\right)$ we have
  \begin{equation}\label{eq:pi_tail}
    \pi(B(0,a)^c) \leq 2d\exp \left(-\frac{a}{2\sqrt{d}(r+1)}\log\frac{a}{2\sqrt{d}M}\right).
  \end{equation}
\end{lemat}

\dowod It follows from \cite[Proposition 11.10]{Sato} that the measure 
$$
  \pi_j(A) = \pi(\{y\in\R:\: y_j\in A\}),\, A\subset\RR,
$$
is infinitely divisible  distribution on $\RR$ with the generating triplet $(0,\nu_0^j,\gamma_j)$, where
$$
  \nu_0^j(A) = \nu_0(\{y\in\R:\: y_j\in A\}),\, A\subset\RR,
$$
and
$$
  \gamma_j = \int y_j \left(\indyk{[-1,1]}(y_j)-\indyk{B(0,1)}(y)\right)\,\nu_0(dy).
$$
Therefore by \cite[Lemma 26.4]{Sato} we obtain
$$
  \pi(\{y\in\R:\:y_j>b\}) \leq \exp\left(-\int_{\xi_0^j}^b \theta_j(\xi)d\xi\right),
$$
for every $b\in (\xi_0^j,\infty)$, where $\theta_j(\xi)$ is the inverse function of
$$
  \psi_j(u) = \int y_j \left(e^{uy_j}-\indyk{B(0,1)}(y)\right)\,\nu_0(dy),\quad u\in\RR.
$$
We note that $\psi_j$ is continuous and increasing,
$$
  \psi_j(0) = \xi_0^j,
$$
and $\lim_{u \to \infty} \psi_j(u) < \infty$ if and only if 
$$
  \nu_0\left(\{y:\:y_j>0\}\right)=0,\, \mbox{ and }\, \int_{-r<y_j<0} |y_j|\,\nu_0(dy)<\infty.
$$
In this case $\supp \pi_j \subset (-\infty,0)$ (see \cite[Corollary 24.8]{Sato}). If
$$
  \lim_{u \to \infty} \psi_j(u) = \infty,
$$
then $\theta_j$ is well-defined on $(\xi_0^j,\infty)$. For every $\xi\in(\xi_0^j,\infty)$ we have
\begin{eqnarray*}
  \xi = \psi_j(\theta_j(\xi))) 
  &   =  & \int y_j \left(e^{\theta_j(\xi)y_j}-\indyk{B(0,1)}(y)\right)\,\nu_0(dy) \\
  &   =  & \int y_j \left(e^{\theta_j(\xi)y_j}-1\right)\,\nu_0(dy) + \xi_0^j \\
  & \leq & e^{r\theta_j(\xi)}\theta_j(\xi)\int y_j^2\, \nu_0(dy) + \xi_0^j \leq \frac{1}{e}e^{\theta_j(\xi)(r+1)}M_j+\xi_0^j.
\end{eqnarray*}
We obtain
$$
  \log(\xi-\xi_0^j)\leq \theta_j(\xi)(r+1)-1+\log M_j,
$$
which yields
$$
  \theta_j(\xi) \geq \frac{1}{r+1}\log\left(e\frac{\xi-\xi_0^j}{M_j}\right).
$$
It follows that
\begin{eqnarray*}
  \pi(\{y\in\R:\:y_j>b\}) 
  & \leq & \exp\left[-\int_{\xi_0^j}^b \frac{1}{r+1}\log\left( e\frac{\xi-\xi_0^j}{M_j}\right)d\xi\right] \\
  & = & \exp\left[-\frac{b-\xi_0^j}{r+1}\log\frac{b-\xi_0^j}{M_j}\right],\quad b>0 \vee \xi_0^j.
\end{eqnarray*}
Similarly,
\begin{eqnarray*}
  \pi(\{y\in\R:\:y_j<-b\})
  &  =   & \pi(\{y\in\R:\:-y_j> b\}) \\
  & \leq & \exp\left[-\int_{-\xi_0^j}^b \frac{1}{r+1}\log\left(e\frac{\xi+\xi_0^j}{M_j}\right)d\xi\right] \\
  & = & \exp\left[-\frac{b+\xi_0^j}{r+1}\log\frac{b+\xi_0^j}{M_j}\right],\quad b> 0\vee (-\xi_0^j).
\end{eqnarray*}
We get
\begin{eqnarray*}
  \pi\left(B(0,a)^c\right)
  & \leq & \sum_{j=1}^d \pi\left(\{y\in\R:\: |y_j| > \frac{a}{\sqrt{d}}\}\right) \\
  & \leq & \sum_{j=1}^d \exp\left[-\frac{\frac{a}{\sqrt{d}}-\xi_0^j}{r+1}\log \frac{\frac{a}{\sqrt{d}}-\xi_0^j}{M_j}\right] \\
  &      & + \sum_{j=1}^d \exp\left[-\frac{\frac{a}{\sqrt{d}}+\xi_0^j}{r+1}\log \frac{\frac{a}{\sqrt{d}}+\xi_0^j}{M_j}\right],
           \quad a>\sqrt{d}\,\xi_0.
\end{eqnarray*}
For $a>2\sqrt{d}\,\xi_0$ we have
$$
  \frac{a}{\sqrt{d}}\pm \xi_0^j \geq \frac{a}{\sqrt{d}} - \xi_0 \geq \frac{a}{\sqrt{d}}-\frac{a}{2\sqrt{d}}=\frac{a}{2\sqrt{d}}.
$$
The function $f(s)=s\log s$ is increasing on $[1/e,\infty)$, which yields
$$
  \exp\left[-\frac{\frac{a}{\sqrt{d}} \pm \xi_0^j}{r+1}\log \frac{\frac{a}{\sqrt{d}} \pm \xi_0^j}{M_j}\right]
  \leq \exp \left(-\frac{a}{2\sqrt{d}(r+1)}\log\frac{a}{2\sqrt{d}M}\right),
$$
provided $a\geq 2\sqrt{d}\left(\xi_0+\frac{M}{e}\right)$, and the lemma follows.
\qed
\\ 

\begin{lemat}\label{lm:pi_estimate} Assume that $\pi$ has a density $p(x)$ and let $m_0$, $m_1$ be such that
  \begin{equation}\label{eq:def_mi}
    \sup_{x\in\R} p(x)\leq m_0,\quad \max_{|\beta|=1}\sup_{x\in\R} |\partial^\beta p(x)| \leq m_1.
  \end{equation}
  Then there exist constants $c_1,c_2,c_3$ such that
  \begin{equation}\label{eg:pi_estimate}
    p(x) \leq c_1 m_1^{d/(d+1)}\exp\left(-\frac{c_2|x|}{r+1}\log\frac{c_3 |x|}{M}\right),
  \end{equation}
  for $|x|>\max\left\{4\sqrt{d}\left(\xi_0+M/e\right),\frac{m_0}{m_1\sqrt{d}}\right\}$.
\end{lemat}
\dowod It follows from Lagrange's theorem that for every pair $x,y\in\R$ there exists $\theta\in [0,1]$ such that 
$$
  p(y)=p(x)+\nabla \scalp{p(x+\theta(y-x))}{(y-x)},
$$
and if $|y-x| \leq \frac{p(x)}{2m_1\sqrt{d}}$ then we obtain
\begin{equation}\label{eq:p_est}
  p(y)\geq p(x) - m_1\sqrt{d}\,|y-x| \geq \frac{1}{2}p(x).
\end{equation}
Let $|x|>\max\left\{4\sqrt{d}\left(\xi_0+M/e\right),\frac{m_0}{m_1\sqrt{d}}\right\}$. By Lemma \ref{lm:pi_tail} 
and (\ref{eq:p_est}) we get
\begin{eqnarray*}
  2d\exp \left(-\frac{|x|/2}{2\sqrt{d}(r+1)}\log\frac{|x|/2}{2\sqrt{d}M}\right)
  & \geq & \pi\left(B(0,|x|/2)^c\right) \\
  & \geq & \pi\left(B(x,\frac{p(x)}{2m_1\sqrt{d}})\right) \\
  &   =  & \int_{B\left(x,\frac{p(x)}{2m_1\sqrt{d}}\right)} p(y)\, dy \\
  & \geq & \frac{1}{2} p(x) c \left(\frac{p(x)}{m_1 \sqrt{d}}\right)^d,
\end{eqnarray*}
and we obtain
$$
  p(x) \leq c_1 m_1^{\frac{d}{d+1}} \exp\left(-\frac{c_2|x|}{r+1}\log\frac{c_3|x|}{M}\right).
$$
\qed

\section{Large jumps}\label{sc:Large_Jumps}

We assume in this section that $\nu$ is a L\'evy measure which satisfies (\ref{eq:nu_intro}), (\ref{eq:q_intro}), (\ref{eq:phi_intro}), 
and (\ref{eq:mu_gamma_measure}). For $r>0$ we denote $\nubounded{r}(dy)= \nu_{r,\infty}(dy)= \indyk{B(0,r)^c}(y)\,\nu(dy)$.

For a set $A\subset\R$ we denote $\delta(A)=\dist(0,A)=\inf\{|y|:\:y\in A\}$ and
$\diam(A)=\sup\{|y-x|:\:x,y\in A\}$. 

For every set $A\subset \R$ such that $\delta(A)>0$ and every $x_0\in A$ we have
\begin{eqnarray}\label{eq:nu_r_est}
  \nubounded{r}(A) 
  & \leq & \int_\sfera \int_{\delta(A)}^\infty \indyk{B(x_0,\diam(A))}(s\theta) 
           s^{-1-\alpha}q(s)\phi(s)\, ds\mu(d\theta) \nonumber \\
  & \leq & \mu\left(B\left(\frac{x_0}{|x_0|},\frac{2\diam(A)}{|x_0|}\right)\right)
           2\diam(A)(\delta(A))^{-1-\alpha}q(\delta(A))\phi(\delta(A)) \nonumber \\
  & \leq & c_1 \left(\frac{2\diam(A)}{|x_0|}\right)^{\gamma-1}(\delta(A))^{-1-\alpha}q(\delta(A))\phi(\delta(A))\diam(A) \nonumber \\
  & \leq & c_2 (\delta(A))^{-\gamma-\alpha}q(\delta(A))\phi(\delta(A)) \left(\diam(A)\right)^\gamma.
\end{eqnarray}

We will extend the estimate to convolutional powers of $\nubounded{r}$. We let 
$$
  \psi(r)=|\mu|\phi(0)\int_r^\infty s^{-1-\alpha}q(s)\frac{\phi(s)}{\phi(s/2)} \, ds,\quad r\in (0,\infty),
$$
and we note that
$$
  |\nubounded{r}| \leq |\mu|\int_r^\infty s^{-1-\alpha}q(s)\phi(s) \, ds \leq \psi(r),\quad r\in (0,\infty).
$$

\begin{lemat}\label{lm:nu_r_n_est} 
There exists a constant $c$ such that 
  \begin{equation}\label{eq:nu_r_n_est}
    \nubounded{r}^{n*}(A) \leq c^n \left(\psi(r)\right)^{n-1} \left(\delta(A)\right)^{-\gamma-\alpha}
    q(\delta(A))\phi(\delta(A)/2) \left(\diam(A)\right)^\gamma.
  \end{equation}
for every set $A$ such that $\delta(A)>0$ and all $n\in\N$, $r>0$.
\end{lemat}

\dowod We use induction. Let (\ref{eq:nu_r_n_est}) hold for some $n\in\N$ and constant $c_0$ and $A$ be a set such that $\delta(A)>0$. For $y\in\R$ we denote $D_y=\{z\in\R:\: |z| > \frac{1}{2}|z+y|\}=B\left(\frac{1}{3}y,\frac{2}{3}|y|\right)^c$. We have
\begin{eqnarray*}
  \nubounded{r}^{(n+1)*}(A) 
  &  =  & \int \nubounded{r}(A-y)\, \nubounded{r}^{n*}(dy) \\
  &  =  & \int \nubounded{r}\left((A-y)\cap D_y\right)\, \nubounded{r}^{n*}(dy) + \int\nubounded{r}\left((A-y)\cap D_y^c\right)\, \nubounded{r}^{n*}(dy) \\
  &  =  & I + II.
\end{eqnarray*}

We note that for $z\in (A-y)\cap D_y$ we have $z+y\in A$ and $|z|>\frac{1}{2} |z+y|$, therefore $|z| >\frac{1}{2} \delta(A)$ and 
$\delta((A-y)\cap D_y)>\frac{1}{2} \delta(A)$. Furthermore, $\diam((A-y)\cap D_y) \leq  \diam(A)$ and using (\ref{eq:nu_r_est}) 
and (\ref{eq:q_intro_2}) we obtain
\begin{eqnarray*}
  I
  & \leq & c_1 2^{\gamma+\alpha+\eta}(\delta(A))^{-\gamma-\alpha}q(\delta(A))
           \phi\left(\delta(A)/2\right) \left(\diam(A)\right)^\gamma |\nubounded{r}^{n*}| \\
  & \leq & c_1 2^{\gamma+\alpha+\eta} \left(\psi(r)\right)^n (\delta(A))^{-\gamma-\alpha}q(\delta(A))\phi\left(\delta(A)/2\right)
           \left(\diam(A)\right)^\gamma.
\end{eqnarray*}

We have
\begin{eqnarray*}
  II
  &  =   & \int\int \indyk{A-y}(z) \indyk{D_y^c}(z)\, \nubounded{r}(dz)\nubounded{r}^{n*}(dy) \\
  &  =   & \int\int \indyk{A-z}(y)\indyk{B(-z,2|z|)^c}(y) \, \nubounded{r}^{n*}(dy)\nubounded{r}(dz)\\
  &  =   & \int \nubounded{r}^{n*}\left((A-z)\cap B(-z,2|z|)^c\right)\,\nubounded{r}(dz),
\end{eqnarray*}

Let $y\in V_z :=(A-z)\cap B(-z,2|z|)^c$. We then have $y+z\in A$, so $|y+z|\geq \delta(A)$,
and $|y+z|\geq 2|z|$. Furthermore $|y|\geq |y+z|-|z|$ and this yields 
$$
  \delta(V_z) \geq \inf_{y\in V_z} |y+z| - |z| \geq \left(\delta(A)\vee 2|z|\right)-|z| \geq \frac{1}{2}\delta(A),
$$
and by (\ref{eq:q_intro}), (\ref{eq:phi_intro}), (\ref{eq:q_intro_2}) and induction we get
\begin{eqnarray*}
  II
  & \leq & \int c_0^n \left(\psi(r)\right)^{n-1} \left(\frac{1}{2}\delta(A)\right)^{-\gamma-\alpha}
           q\left(\frac{1}{2}\delta(A)\right) \times \\
  &      & \times \,\phi\left(\frac{1}{2} \left(\delta(A)\vee 2|z|\right)-\frac{1}{2} |z|\right) \left(\diam(A)\right)^\gamma\, \nubounded{r}(dz) \\
  & \leq & c_0^n c_2 \left(\psi(r)\right)^{n-1} 2^{\gamma+\alpha+\eta}\left(\delta(A)\right)^{-\gamma-\alpha}
           q\left(\delta(A)\right) \left(\diam(A)\right)^\gamma\ \times\\
  &      & \times \kappa_2 \int \frac{\phi\left(\delta(A)/2\right)}{\phi(|z|/2)}\, \nubounded{r}(dz) \\
  & \leq & c_0^n c_2 \left(\psi(r)\right)^{n-1} 2^{\gamma+\alpha+\eta}\left(\delta(A)\right)^{-\gamma-\alpha}
           q\left(\delta(A)\right) \phi\left( \delta(A)/2\right) \left(\diam(A)\right)^\gamma \times \\
  &      & \times \kappa_2 |\mu| \int_r^\infty \frac{1}{\phi(s/2)} s^{-1-\alpha}q(s)\phi(s)\, ds \\
  & =    & c_0^n \kappa_2 c_2 2^{\gamma+\alpha+\eta} \left(\psi(r)\right)^{n-1} \left(\delta(A)\right)^{-\gamma-\alpha}
           q\left(\delta(A)\right) \phi\left( \delta(A)/2\right) \left(\diam(A)\right)^\gamma \frac{\psi(r)}{\phi(0)}.
\end{eqnarray*}
The lemma follows by taking $c_0 \geq \max\{1,2^{\gamma+\alpha+\eta}(c_1+\kappa_2 c_2/\phi(0))\}$.

\qed
\\
\begin{wniosek} There exists a constant $c$ such that for every $x\neq 0$, $n\in\N$ and $\rho<\frac{1}{2}|x|$ we have
  \begin{equation}\label{eq:nu_dla_kul}
    \nubounded{r}^{n*}(B(x,\rho)) \leq c^n \left(\psi(r)\right)^{n-1} |x|^{-\gamma-\alpha}q(|x|)
    \phi\left(|x|/4 \right) \rho^\gamma.
  \end{equation}
\end{wniosek}
We consider the probability measures $\{\bar{P}^r_t,\; t\geq 0\}$ such that
\begin{equation}\label{eq:FPbar}
  \Fourier(\bar{P}^r_t)(\xi) =
  \exp\left(t \int (e^{i\scalp{\xi}{y}}-1)\,
  \nubounded{r}(dy)\right)\, ,
  \quad \xi\in\R\, .
\end{equation}
Note that
\begin{eqnarray}\label{eq:exp}
  \bar{P}^r_t=P^{r,\infty}_t* \delta_{-ta_{r,1}} 
  &  =  & \exp(t(\nubounded{r}-|\nubounded{r}|\delta_0)) =  \sum_{n=0}^\infty \frac{t^n\left(\nubounded{r}-|\nubounded{r}|\delta_0)\right)^{n*}}{n!} \\
  &  =  & e^{-t|\nubounded{r}|} \sum_{n=0}^\infty \frac{t^n\nubounded{r}^{n*}}{n!}\,,\quad t\geq 0\, . \nonumber
\end{eqnarray}
We obtain
\begin{wniosek}\label{cor:sigma_estimate} There is $c$ such that
  \begin{equation}\label{eq:sigma_estimate}
    \bar{P}^r_t \left(B(x,\rho)\right) \leq c t e^{t(c\psi(r)-|\nubounded{r}|)}  |x|^{-\gamma-\alpha}q(|x|)
    \phi\left(|x|/4\right) \rho^\gamma, 
  \end{equation}
  for every $x\in\R\setminus\{0\}$, and $\rho<\frac{1}{2} |x|$.
\end{wniosek}

\section{Proof of Theorem \ref{th:main}}\label{sc:Main_Theorem}

We begin with the following technical lemma.

\begin{lemat}\label{lm:tech} Let $a>v\geq 0$. If $f:\:(0,\infty) \mapsto \RR$ is a nonnegative function such that
  \begin{equation}\label{eq:tech_ass}
    \int_0^r s^a f(s)\, ds \leq c_1 r^{v},\quad r\geq r_0,
  \end{equation}
  for some constants $c_1,r_0>0$, then
  $$
    \int_r^\infty f(s)\, ds \leq \frac{c_1 a}{a-v} r^{v-a},\quad r\geq r_0.
  $$
\end{lemat}

\dowod For $r>r_0$ by (\ref{eq:tech_ass}) we have
$$
  \int_r^\infty \frac{1}{t^{a+1}}   \left( \int_0^t s^a f(s)\, ds\right) dt \leq \int_r^\infty c_1 t^{v-a-1} dt = \frac{c_1}{a-v} r^{v-a}.
$$
Furthermore, changing the order of integration we obtain
\begin{eqnarray*}
  \int_r^\infty \frac{1}{t^{a+1}}   \left( \int_0^t s^a f(s)\, ds\right) dt 
  &  =   & \int_0^r s^a f(s) \int_r^\infty \frac{1}{t^{a+1}} \, dt \, ds+ \int_r^\infty s^a f(s) \int_s^\infty \frac{1}{t^{a+1}} \, dt\, ds \\
  &  =   & \frac{1}{ar^a} \int_0^r s^a f(s) \, ds + \frac{1}{a} \int_r^\infty f(s)\, ds \\
  & \geq & \frac{1}{a} \int_r^\infty f(s)\, ds,
\end{eqnarray*}
and the lemma follows.

\qed
\\

For $r>0$ we denote $\tilde{\nu}_r(dy)=\nu_{0,r}(dy)=\indyk{B(0,r)}(y)\nu(dy)$. 
We consider the semigroup of measures $\{\tilde{P}^r_t,\; t\geq 0\}$ such that
$$
  \Fourier(\tilde{P}^r_t)(\xi) 
    =     \exp\left(t \int \left(e^{i\scalp{\xi}{y}}-1-i\scalp{\xi}{y}\right)
            \tilde{\nu}_r(dy)\right)\, ,  \quad \xi\in\R\,.
$$
By (\ref{eq:Ass0}) we get
\begin{eqnarray}\label{eq:FTildePEstimate}
  |\Fourier(\tilde{P}^r_t)(\xi)| 
  &   =  & \exp\left(-t\int_{|y|<r}
           (1-\cos(\scalp{y}{\xi}))\,
           \nu(dy)\right) \nonumber \\
  &   =  & \exp\left(-t\left(\Re(\Phi(\xi))-\int_{|y|\geq r}
           (1-\cos(\scalp{y}{\xi}))\,
           \nu(dy)\right)\right) \nonumber \\
  & \leq & \exp(-t\Re(\Phi(\xi)))\exp(2t\nu(B(0,r)^c)) \nonumber \\
  & \leq & \exp\left(-ct\left(|\xi|^\alpha\wedge |\xi|^\beta\right)\right)\exp\left(2t\nu(B(0,r)^c)\right),
           \quad \xi\in\R.
\end{eqnarray}
It follows  that for every $r>0$ and $t>0$ the measure $\tilde{P}^r_t$ is absolutely continuous with respect to
the Lebesgue measure with smooth density, say, $\tilde{p}^r_t$.

Let $h(t)=t^{1/\alpha}\wedge {t^{1/\beta}}$. Informally, $h(t)$ gives a correct (time-dependent) threshold for truncating the L\'evy measure of convolutional exponents considered as approximations of the original transition probability $p_t$. We will often use $\tilde{P}^r_t$ and $\tilde{p}^r_t$ with $r=h(t)$ and
for simplification we will denote 
$$ 
  \tilde{P}_t=\tilde{P}^{h(t)}_t\,\, \mbox{and}\,\, \tilde{p}_t=\tilde{p}^{h(t)}_t.
$$

\begin{lemat}\label{lm:1} If (\ref{eq:nu_intro}), (\ref{eq:beta_intro}) and (\ref{eq:Ass0}) hold then
  \begin{equation}\label{eq:p_tilde_estimate}
    \tilde{p}_t(x) \leq c_1 [h(t)]^{-d} \exp\left( -c_2\frac{|x|}{h(t)} \log\left( 1+ c_3\frac{|x|}{h(t)}\right)\right),
  \end{equation}
  for every $x\in\R\setminus\{0\}$, $t>0$.
\end{lemat}

\dowod Let $g_t(y)=(h(t))^d\tilde{p}_t\left(h(t)y\right)$. It follows from (\ref{eq:nu_intro}) that for $r<1$ we have 
$\nu(B(0,r)^c)\leq \left( |\mu|q(0)\phi(0)/\alpha\right)r^{-\alpha}$, and from (\ref{eq:beta_intro}) and Lemma \ref{lm:tech} with $a=2$ 
and $v=2-\beta$ for $r>1$ we get 
\begin{eqnarray*}
  \nu(B(0,r)^c) 
  & \leq &  |\mu| \int_r^\infty s^{-1-\alpha}q(s)\phi(s)\, ds \\
  & \leq & \phi(0)|\mu| \int_r^\infty s^{-1-\alpha}q(s)\frac{\phi(s)}{\phi(s/2)}\,ds \leq c_1 r^{-\beta}.
\end{eqnarray*}
From (\ref{eq:FTildePEstimate}) 
for every $j\in\{1,\dots,d\}$ we obtain
\begin{eqnarray*}
  \left|\frac{\partial g_t}{\partial y_j}(y)\right|
  &   =  & (h(t))^{d+1} \left|(2\pi)^{-d}\int_{\R}
           (-i)\xi_j e^{-i\scalp{h(t)y}{\xi}}
           \Fourier(\tilde{p}_t)(\xi)d\xi\right|  \\
  & \leq & c_2 (h(t))^{d+1} \left(\int_{|\xi|\leq 1}|\xi_j|e^{-c_3 t |\xi|^\beta}\,d\xi + 
                     \int_{|\xi|>1}|\xi_j|e^{-c_3 t|\xi|^\alpha}\,d\xi \right)  \\
  &   =  & c_2 (h(t))^{d+1} \left( t^{\frac{-1-d}{\beta}} \int_{|u|\leq t^{1/\beta}}|u_j|e^{-c_3 |u|^\beta}\,du+
                      t^{\frac{-1-d}{\alpha}} \int_{|u| > t^{1/\alpha}}|u_j|e^{-c_3 |u|^\alpha}\,du \right)   \\
  & \leq & c_4 (h(t))^{d+1} \left( t^{\frac{-1-d}{\beta}} + t^{\frac{-1-d}{\alpha}}\right)  \\
  & \leq & 2c_4.
\end{eqnarray*}
Similarly we get $g_t(y)\leq c_5$ for all $y\in\R$, $t>0$.

We consider the infinitely divisible distribution $\pi_t(dy)=g_t(y)\, dy$. We note that
\begin{eqnarray*}
  \Fourier(\pi_t)(\xi) 
  &  =   &  \exp\left(t \int \left(e^{i\scalp{\xi(h(t))^{-1}}{y}}-1-i\scalp{\xi(h(t))^{-1}}{y}\indyk{B(0,h(t))}(y)\right)
            \tilde{\nu}_{h(t)}(dy)\right)\\
  &  =   & \exp\left( \int \left(e^{i\scalp{\xi}{y}}-1-i\scalp{\xi}{y}\indyk{B(0,1)}(y)\right)
            \lambda_t (dy)\right)\, ,  \quad \xi\in\R\, ,
\end{eqnarray*}
where $\lambda_t(A)=t\tilde{\nu}_{h(t)}(h(t)A)$ is the L\'evy measure of $\pi_t$. 
By (\ref{eq:beta_intro}) we have
\begin{eqnarray*}
  \int |y|^2\, \lambda_t(dy) 
  &   =  &  t \int (|y|/h(t))^2\,\tilde{\nu}_{h(t)}(dy) \\
  & \leq &  |\mu| t(h(t))^{-2} \int_0^{h(t)} s^{1-\alpha}q(s)\phi(s)\, ds \\
  & \leq &   c_6,
\end{eqnarray*}
and
\begin{eqnarray*}
  \int_{|y|>1} y_j \lambda_t(dy) 
  &  =   &  t (h(t))^{-1}\int_{B(0,h(t))^c} y_j \tilde{\nu}_{h(t)}(dy) = 0. 
\end{eqnarray*}
It follows from Lemma \ref{lm:pi_estimate} that
$$
  g_t(y) \leq c_7 \exp \left( -c_8|y| \log\left( c_9|y|\right)\right) \leq c_{10} \exp \left( -c_8|y| \log\left(1+ c_9|y|\right)\right),
$$
for $y\in\R\setminus\{0\}$, and this yields
$$
  \tilde{p}_t(x) \leq c_{10} (h(t))^{-d} \exp\left( -\frac{c_8|x|}{h(t)} \log\left(1+ \frac{c_9|x|}{h(t)}\right)\right),
$$
for $x\in\R\setminus\{0\}$, $t>0$.
\qed

We will prove now the on-diagonal estimates of the density $p_t$. We observe that Nash-type inequalities (\cite{CKS}) can be used to prove such estimates in the case of (symmetric) Markov semigroups. The present situation of convolutional semigroups is, however, simpler, and we can use Fourier transform instead.

We note that the existence (and smoothness) of the densities follows from (\ref{eq:Ass0}) 
because (\ref{eq:Ass0}) yields that the Fourier transform of $P_t$
decays faster in infinity than every negative power of $|\xi|$.

\begin{lemat}\label{lm:pBasicEstimate} If (\ref{eq:Ass0}) holds then there exists $c$ such that
  $$
    p_t(x) \leq c [h(t)]^{-d},\quad x\in\R,\, t > 0.
  $$
\end{lemat}
\dowod 
By (\ref{eq:Ass0}) we have
\begin{eqnarray*}
	p_t(x) &  =   & (2\pi)^{-d} \int_{\R} e^{-i\scalp{x}{\xi}}e^{-t\Phi(\xi)}\, d\xi \\
	       & \leq & (2\pi)^{-d} \int_{\R} |e^{-t\Phi(\xi)}| \,d\xi = (2\pi)^{-d} \int_{\R} e^{-t\Re(\Phi(\xi))} \,d\xi \\
	       & \leq & (2\pi)^{-d} \left( \int_{|\xi|\leq 1} e^{-c_1 t|\xi|^\beta}\,d\xi 
	                                  + \int_{|\xi|\geq 1} e^{-c_1 t|\xi|^\alpha}\,d\xi \right) \\
	       &   =  & (2\pi)^{-d} \left( t^{-d/\beta} \int_{|u| \leq t^{1/\beta}} e^{-c_1|u|^\beta}\,du
                    + t^{-d/\alpha} \int_{|u|\geq t^{1/\alpha}}e^{-c_1|u|^\alpha}\,du\right)\\
         & \leq & c_2 \left(  t^{-d/\beta}+ t^{-d/\alpha}\right)\\
         & \leq & 2c_2 [h(t)]^{-d} ,\quad x\in\R,\, t>0.
\end{eqnarray*}
\qed

\dowod {\it of Theorem \ref{th:main}}. We have
\begin{displaymath}
  P_t=\tilde{P}^r_t \ast \bar{P}^r_t \ast \delta_{t b_r}\,,\quad t\geq 0,
\end{displaymath}
where $\bar{P}^r_t$ is defined by (\ref{eq:FPbar}) and $b_r$ by (\ref{eq:def_br}).
Of course
\begin{equation}\label{eq:wpt}
  p_t= \tilde{p}^r_t * \bar{P}^r_t\ast \delta_{t b_r} \,, \quad t>0.
\end{equation}

We will denote 
$$
  \bar{P}_t=\bar{P}_t^{h(t)}.
$$
By (\ref{eq:beta_intro}) and Lemma \ref{lm:tech} we have
\begin{eqnarray*}
  \psi(h(t)) 
  &  =   & |\mu|\phi(0)\int_{h(t)}^\infty s^{-1-\alpha}q(s)\frac{\phi(s)}{\phi(s/2)} \, ds \\
  & \leq & \frac{c}{t},\quad t>0.
\end{eqnarray*}
It follows from Corollary \ref{cor:sigma_estimate} that
\begin{equation}\label{eq:bar_P_schonwieder}
  \bar{P}_t(B(x,\rho)) \leq c t |x|^{-\gamma-\alpha}q(|x|)
    \phi\left(|x|/4\right) \rho^\gamma,
\end{equation}
for $\rho\leq \frac{1}{2}|x|$ and $t>0$.

We denote
$$
  g(s)=e^{-c_2s\log(1+c_3s)},\quad s\geq 0,
$$
where constants $c_2,c_3$ are given by (\ref{eq:p_tilde_estimate}). We note that $g$ is decreasing, continuous on $[0,\infty)$ and $g(s)\leq c s^{-2\gamma}$, for some
$c>0$, which yields that the inverse function $g^{-1}:\: (0,1]\to [0,\infty)$ exists, is decreasing, and $g^{-1}(s)\leq \left(c/s\right)^{1/(2\gamma)}$.
In particular
$$
  \int_0^1 \left(g^{-1}(s)\right)^{\gamma}\,ds < \infty.
$$
Using Lemma \ref{lm:1} and (\ref{eq:bar_P_schonwieder}) we obtain
\begin{eqnarray*}
  \tilde{p}_t * \bar{P}_t (x)
  &  =   & \int \tilde{p}_t(x-y)\, \bar{P}_t(dy) \\
  & \leq & \int c [h(t)]^{-d} g(|x-y|/h(t)) \, \bar{P}_t(dy) \\
  &   =  & c [h(t)]^{-d} \int \int_0^{g(|x-y|/h(t))} \, ds\, \bar{P}_t(dy) \\
  &   =  & c [h(t)]^{-d} \int_0^1 \int \indyk{\{y\in\R:\: g(|x-y|/h(t))>s\}} \, \bar{P}_t(dy) ds \\
  &   =  & c [h(t)]^{-d} \int_0^1 \bar{P}_t\left(B(x,h(t)g^{-1}(s))\right) ds \\
  & \leq & c [h(t)]^{-d} \left(\int_{g(\frac{|x|}{2h(t)})}^1 t |x|^{-\gamma-\alpha}q(|x|)\phi\left(|x|/4\right) 
           \left(h(t)g^{-1}(s)\right)^\gamma\, ds + \int_0^{g(\frac{|x|}{2h(t)})}\, ds\right) \\
  & \leq & c [h(t)]^{-d} \left( t[h(t)]^{\gamma} |x|^{-\gamma-\alpha}q(|x|)\phi\left(|x|/4\right) \int_0^1  
           \left(g^{-1}(s)\right)^\gamma\, ds + g\left(\frac{|x|}{2h(t)}\right) \right) \\
  &  =   & c [h(t)]^{-d} \left( t[h(t)]^{\gamma} |x|^{-\gamma-\alpha}q(|x|)\phi\left(|x|/4\right) 
           + g\left(\frac{|x|}{2h(t)}\right) \right). \\
\end{eqnarray*}

This and Lemma \ref{lm:pBasicEstimate} yield
\begin{eqnarray*}
  p_t(x+t b_{h(t)}) 
  &  =   & \int \tilde{p}_t * \bar{P}_t (x+t b_{h(t)}-y) \delta_{t b_{h(t)}}(dy) \\
  &  =   & \tilde{p}_t * \bar{P}_t (x) \\
  & \leq & c [h(t)]^{-d} \min\left\{1, t[h(t)]^{\gamma}
          |x|^{-\gamma-\alpha}q(|x|)\phi\left(|x|/4\right) + \,g\left(\frac{|x|}{2h(t)}\right) \right\}, \\
  & \leq & c [h(t)]^{-d}\left( \min\left\{1, t[h(t)]^{\gamma}
          |x|^{-\gamma-\alpha}q(|x|)\phi\left(|x|/4\right)\right\} + \,g\left(\frac{|x|}{2h(t)}\right) \right).
\end{eqnarray*}
For $t>1$ we have
\begin{eqnarray*}
  p_t(x+t b_{t^{1/\beta}}) 
  & \leq & c t^{-d/\beta}\left( \min\left\{1, t^{1+\frac{\gamma}{\beta}}
          |x|^{-\gamma-\alpha}q(|x|)\phi\left(|x|/4\right)  \right\}
           + g\left(\frac{|x|}{2t^{1/\beta}}\right)\right),\quad x\in\R.
\end{eqnarray*}
We observe that the function $1/\phi$ is submultiplicative and locally bounded therefore by \cite[Lemma 25.5]{Sato} there exist $c_1, c_2$ such that $\phi(s)\geq c_1 e^{-c_2 s}$ and this and (\ref{eq:q_intro_2}) yield
$$
  g(s)\leq c s^{-\gamma-\alpha}q(s)\phi\left(s\right),\quad s>0.
$$
For $t\in (0,1)$ we get
\begin{eqnarray*}
  p_t(x+t b_{t^{1/\alpha}}) 
  & \leq & c t^{-d/\alpha} \min\left\{1, t^{1+\frac{\gamma}{\alpha}}
          |x|^{-\gamma-\alpha}q(|x|)\phi\left(|x|/4\right) 
            \right\},\quad x\in\R.
\end{eqnarray*}
\qed

\section{Examples}\label{sec:examples}
\subsection{Stable processes}

We call a measure $\lambda$ {\it degenerate} if there is a proper
linear subspace $M$ of $\R$ such that $\supp(\lambda)\subset M$; 
otherwise we call $\lambda$ {\it nondegenerate}. 

Let $\mu$ be a nondegenerate bounded measure on $\sfera$, $\alpha\in(0,2)$,  and
\begin{equation}\label{eq:L_stable}
  \nu(A)=\int_0^\infty \int_{\sfera} \indyk{A}(s\theta)s^{-1-\alpha}\,ds \mu(d\theta),\quad A\subset\R.
\end{equation}
Let
$$
  b = \left\{
  \begin{array}{ccc}
    \int_{|y|<1} y \, \nu(dy)   & \mbox{  if  } & \alpha < 1,\\
    0                           & \mbox{  if  } & \alpha=1, \\
    - \int_{1<|y|} y \, \nu(dy) & \mbox{  if  } & \alpha > 1.
  \end{array}\right.
$$
If $\alpha=1$ then we assume additionally that
$$
  \int_{\sfera} \theta\,\mu(d\theta) = 0.
$$
The corresponding semigroup $\{P_t:\:t\geq 0\}$ is called the $\alpha$--stable semigroup and the stochastic process
$\{X_t:\: t\geq 0\}$ is the $\alpha$--stable L\'evy process. We have (see, e.g., \cite[Theorem 14.10]{Sato})
\begin{equation}\label{eq:Phi_stable_0}
  \Phi(\xi) = a_\alpha \int_\sfera |\scalp{\xi}{\theta}|^{\alpha}\left(1-i\tan \frac{\pi\alpha}{2}\sgn \left(\scalp{\xi}{\theta}\right)\right)\,\mu(d\theta),\quad 
  \mbox{if}\quad \alpha\neq 1,
\end{equation}
and
\begin{equation}\label{eq:Phi_stable_1}
  \Phi(\xi) = a_1 \int_\sfera |\scalp{\xi}{\theta}|
  \left(1+i\frac{2}{\pi} \sgn\left(\scalp{\xi}{\theta}\right)\log|\scalp{\xi}{\theta}|\right)\,\mu(d\theta),\quad 
  \mbox{if}\quad \alpha =1,
\end{equation}
where $a_\alpha=\frac{\pi}{2\sin\frac{\pi\alpha}{2}\Gamma(1+\alpha)}$. From the fact that $\mu$ is nondegenerate it follows easily that $\Re(\Phi(\xi))\geq c |\xi|^\alpha$. Of course $\Phi(a\xi)=a^{\alpha}\Phi(\xi)$ for every $a>0$ and the semigroup has a scaling property, i.e.,
$$
  P_t(A) = P_1(t^{-1/\alpha}A),\quad A\subset\R,
$$
and
$$
  p_t(x) = t^{-d/\alpha}p_1(t^{-1/\alpha}x),\quad t>0,\,x\in\R.
$$
Using Theorem \ref{th:main} for $t=1$ with $q(s)\equiv\phi(s)\equiv 1$, and $\beta=\alpha$ we obtain
\begin{wniosek}
  If $\nu$ is given by (\ref{eq:L_stable}), $\mu$ satisfies (\ref{eq:mu_gamma_measure}) for some $\gamma\in[1,d]$, 
  and $\{P_t:\:t\geq 0\}$ is the stable semigroup determinated 
  by the characteristic exponent $\Phi$ given in (\ref{eq:Phi_stable_0}) and (\ref{eq:Phi_stable_1}) then we have
  $$
    p_t(x) \leq c t^{-d/\alpha} \min\{1,t^{1+\frac{\gamma}{\alpha}}|x|^{-\gamma-\alpha}\},\quad x\in\R, t>0.
  $$
\end{wniosek}

We note that such estimates were first obtained in \cite{BS2007} for symmetric stable semigroups, and in \cite{W} for non-symmetric stable semigroups. Our methods are different than those of \cite{W}.

\subsection{Layered stable processes}
Let $\mu$ be a nondegenerate bounded measure on $\sfera$, and
\begin{equation}\label{eq:L_lay_stable}
  \nu(A)=\int_0^\infty \int_{\sfera} \indyk{A}(s\theta)Q(\theta,s)\,ds \mu(d\theta),\quad A\subset\R,
\end{equation}
where
$$
  Q(\theta,s) \approx s^{-\alpha-1}\indyk{(0,1]}(s)+s^{-m-1}\indyk{(1,\infty)}(s),
$$
and $m > 2$. Let $b\in\R$ and
\begin{eqnarray}\label{eq:Phi_lay}
  \Phi(\xi) 
  &   =  &  - \int \left(e^{i\scalp{\xi}{y}}-1-i\scalp{\xi}{y}\indyk{B(0,1)}(y)\right)\nu(dy) - i\scalp{\xi}{b}\nonumber \\
\end{eqnarray}

A semigroup with the characteristic exponent (\ref{eq:Phi_lay}) is called {\it layered stable semigroup} and was introduced
by C. Houdr\'e and R. Kawai in \cite{HoKawa}.

\begin{wniosek}
  If $\nu$ is given by (\ref{eq:L_lay_stable}), $\mu$ satisfies (\ref{eq:mu_gamma_measure}) for some $\gamma\in[1,d]$, 
  and $\{P_t:\:t\geq 0\}$ is the layered stable semigroup determinated 
  by the characteristic exponent $\Phi$ given in (\ref{eq:Phi_lay}) then the density $p_t$ of $P_t$ satisfies
  \begin{eqnarray*}
  p_t(x+tb_{t^{1/\alpha}}) 
  & \leq & c t^{-d/\alpha} \min\left\{1, t^{1+\frac{\gamma}{\alpha}}
          |x|^{-\gamma-\alpha}(1+|x|)^{\alpha-m}\right\},
\end{eqnarray*}
for $t\in(0,1],\,x\in\R,$
and
\begin{eqnarray*}
  p_t(x+tb_{\sqrt{t}}) 
  & \leq & c_1 t^{-d/2} \left(\min\left\{1, t^{1+\frac{\gamma}{2}}
          |x|^{-\gamma-\alpha}(1+|x|)^{\alpha-m}\right\} 
           + e^{\frac{-c_2|x|}{\sqrt{t}}\log\left(1+\frac{c_3|x|}{\sqrt{t}}\right)} \right),
\end{eqnarray*}
for $t>1$, $x\in\R$.
\end{wniosek}
\dowod It follows from Lemma 10 in \cite{Sztonyk2010_2} that in this case
$$
  \Re\left(\Phi(\xi)\right) \approx |\xi|^2 \wedge |\xi|^\alpha,\quad \xi\in\R.
$$
We let $q(s)=(1+s)^{\alpha-m}$, $\phi\equiv 1$, and $\beta=2$ in Theorem \ref{th:main}.
Of course
$$
  \int_0^r s^{1-\alpha}(1+s)^{\alpha-m}\, ds \leq \int_0^\infty s^{1-\alpha}(1+s)^{\alpha-m}\, ds < \infty,
$$
and (\ref{eq:beta_intro}) is also satisfied with $\beta=2$
and the Corollary follows from Theorem \ref{th:main}.

\qed

\subsection{Tempered stable processes} Let $\mu$ be a as above nondegenerate bounded measure on $\sfera$, and
\begin{equation}\label{eq:L_temp_stable}
  \nu(A)=\int_0^\infty \int_{\sfera} \indyk{A}(s\theta)s^{-1-\alpha} Q(\theta,s)\,ds \mu(d\theta),\quad A\subset\R.
\end{equation}
If $Q(\theta,\cdot)$ is completely monotone for every $\theta\in\sfera$ then a process $\proces$ 
which has the L\'evy measure $\nu$ is calling
{\it tempered stable process}. Basic properties of tempered stable processes were investigated in \cite{Ros07}. We will focus
on functions $Q$ which decay exponentially: we will assume that
$$
  Q(\theta,s) \approx e^{-c s},
$$
for some constant $c>0$.
Let $b\in\R$,
and
\begin{equation}\label{eq:Phi_temp}
  \Phi(\xi) 
     =    - \int \left(e^{i\scalp{\xi}{y}}-1-i\scalp{\xi}{y}\indyk{B(0,1)}(y)\right)\nu(dy) - i\scalp{\xi}{b}.
\end{equation}
Taking $q\equiv 1$ and $\phi(s)=e^{-cs}$ in Theorem \ref{th:main} we obtain the following estimates.
\begin{wniosek}
  If $\nu$ is given by (\ref{eq:L_temp_stable}), $\mu$ satisfies (\ref{eq:mu_gamma_measure}) for some $\gamma\in[1,d]$, 
  and $\{P_t:\:t\geq 0\}$ is the tempered stable semigroup determinated 
  by the characteristic exponent $\Phi$ given by (\ref{eq:Phi_temp}), then the density $p_t$ of $P_t$ satisfies
  \begin{eqnarray*}
  p_t(x+tb_{t^{1/\alpha}}) 
  & \leq & c_1 t^{-d/\alpha} \min\left\{1, t^{1+\frac{\gamma}{\alpha}}
          |x|^{-\gamma-\alpha}e^{-c_2 |x|}\right\},
\end{eqnarray*}
for $t\in(0,1],\,x\in\R,$
and
\begin{eqnarray*}
  p_t(x+tb_{\sqrt{t}}) 
  & \leq & c_3 t^{-d/2} \left(\min\left\{1, t^{1+\frac{\gamma}{2}}
          |x|^{-\gamma-\alpha}e^{-c_4|x|}\right\} 
           + e^{\frac{-c_5|x|}{\sqrt{t}}\log\left(1+\frac{c_6|x|}{\sqrt{t}}\right)} \right),
\end{eqnarray*}
for $t>1$, $x\in\R$.
\end{wniosek}

\subsection{Absolutely continuous and symmetric L\'evy measure}
Let $\nu$ be a symmetric measure, i.e., $\nu(A)=\nu(-A)$, and
\begin{equation}\label{eq:ac_L}
  \nu(dx) \approx |x|^{-d-\alpha} q(|x|)\phi(|x|)\,dx,\quad x\in\R\setminus \{0\},
\end{equation}
where $q$ and $\phi$ satisfy (\ref{eq:q_intro}), (\ref{eq:phi_intro}) and (\ref{eq:beta_intro}) with $\beta=2$.
It follows from Theorem \ref{th:main} here and Theorems 2 and 4 in \cite{Sztonyk2010_2} that the density $p_t$ in this case satisfies
the following estimates
$$
 c_1  \min\left\{t^{-d/\alpha}, t
          |x|^{-d-\alpha}q(|x|)\phi\left(2|x|\right) 
            \right\} \leq
$$
$$
  p_t(x) \leq c_2  \min\left\{t^{-d/\alpha}, t
          |x|^{-d-\alpha}q(|x|)\phi\left(|x|/4\right) 
            \right\},
$$
for $t\in (0,1]$, $x\in\R$  and
$$
  c_3  \min\left\{t^{-d/2}, t
          |x|^{-d-\alpha}q(|x|)\phi\left(2|x|\right)\right\} 
        \leq 
$$
$$
p_t(x) \leq c_4  \left(\min\left\{t^{-d/2}, t
          |x|^{-d-\alpha}q(|x|)\phi\left(|x|/4\right)\right\} 
           + e^{\frac{-c_5|x|}{\sqrt{t}}\log\left(1+\frac{c_6|x|}{\sqrt{t}}\right)} \right),
$$
for $t>1$, $x\in\R$.

In particular for the relativistic $\alpha$--stable L\'evy process, which is investigated, e.g., in \cite{Ryz,KulczSiu,GrzRyz08}, 
and \cite{BMR}, we have
\begin{eqnarray*}
  \nu(D) & =  & c_1 \int_{D} |y|^{-d-\alpha} K_{d,\alpha}(|y|)\,dy \\
         & =  & c_2 \int_{\sfera} \int_0^\infty  \indyk{D}(s\theta) s^{-1-\alpha} K_{d,\alpha}(s)\, ds \sigma(d\theta),
        \quad D\subset\R,
\end{eqnarray*}
where $\sigma$ is the standard isotropic surface measure on $\sfera$ and 
$$
  K_{d,\alpha}(s)=s^{d+\alpha}\int_0^\infty e^{-u}e^{-\frac{s^2}{4u}}u^{\frac{-2-d-\alpha}{2}}\,du,\quad s>0.
$$
We have (see \cite{GrzRyz08})
$$
 K_{d,\alpha}(s) \approx (1+s)^{\frac{d+\alpha-1}{2}}e^{-s}
$$
and in this case we get
$$
 c_1  \min\left\{t^{-d/\alpha}, t
          |x|^{-d-\alpha}e^{-2|x|}
            \right\} \leq
  p_t(x) \leq c_2  \min\left\{t^{-d/\alpha}, t
          |x|^{-d-\alpha}e^{-|x|/5} 
            \right\},
$$
for $t\in (0,1]$, $x\in\R$  and
$$
  c_3  \min\left\{t^{-d/2}, t
          |x|^{-d-\alpha}e^{-2|x|}\right\} 
        \leq 
$$
$$
p_t(x) \leq c_4  \left(\min\left\{t^{-d/2}, t
          |x|^{-d-\alpha}e^{-|x|/5}\right\} 
           + e^{\frac{-c_5|x|}{\sqrt{t}}\log\left(1+\frac{c_6|x|}{\sqrt{t}}\right)} \right),
$$
for $t>1$, $x\in\R$. The sharp estimates of the transition densities of the relativistic process for small $t$ are given
in \cite{ChKum08}.




\end{document}